\documentclass[11pt, twoside]{article}
\usepackage{amsfonts}
\usepackage{amsmath,amssymb}
\usepackage{multicol}
\usepackage{color}
\usepackage{mathtools}
\usepackage{mathtools}
\usepackage{cite}
\newtheorem{theorem}{Theorem}[section]
\newtheorem{lemma}[theorem]{Lemma}

\newtheorem{definition}[theorem]{Definition}
\newtheorem{remark}[theorem]{Remark}

\numberwithin{equation}{section}
\newenvironment{proof}[1][Proof]{\noindent\textbf{#1. }}{\hfill $\Box$}
\allowdisplaybreaks

\makeatletter
\begin{document}
\author{Quanguo Zhang \thanks{Corresponding author:
zhangqg07@163.com.}
\\
{\small Department of Mathematics, Luoyang Normal University,
 Luoyang, Henan 471022, P.R. China}}
\title{\textbf{\Large  Nonexistence of global weak solutions for a wave equation with
nonlinear memory and damping terms
  } } \maketitle
\begin{abstract} In this paper, we study the nonexistence of global weak solutions  for a wave equation with nonlinear memory and damping terms.
We give an answer to an open problem posed in [M. D'Abbicco, A wave equation with structural damping
and nonlinear memory, Nonlinear Differ. Equ. Appl. 21 (2014), 751-773]. Moreover, comparing with the existing results, our results do not require any positivity condition of the initial values. The proof of our results is based on the asymptotic properties of solutions for an integral inequality.


\end{abstract}

\section{Introduction}

\noindent

This paper is mainly concerned with nonexistence of global solutions for the following wave equation with nonlinear memory and damping terms:
\begin{equation}\label{20.1}
\left\{\begin{array}{l}u_{tt}+(-\triangle)^\sigma u+\mu(-\triangle)^\eta u_t={}_0^{}I_t^{1-\gamma}(|u|^{p}),\ \
(t,x)\in (0,T)\times\mathbb{R}^N,\\
u(0,x)=u_0(x),\ \  u_t(0,x)=u_1(x),\ \  x\in \mathbb{R}^N,
\end{array}\right.
\end{equation}
where $0<\sigma\leq 1$, $0\leq \eta\leq 1$, $\mu\in\mathbb{R}$, $\gamma<1$, $p>1$ and ${}_0^{}I_t^{1-\gamma}$ denotes the Riemann-Liouville fractional integral of order $1-\gamma$.

Recently, the heat equation with nonlinear memory term and the wave equation with nonlinear memory term have received extensive attentions. There are
a lot of papers on the global existence and blow-up of solutions for these questions, see \cite{T.Cazenave,chen,chen2, a1,a2,fino,dao,Berbiche2,Andrade} and the
references therein. For example, in \cite{T.Cazenave}, Cazenave et al. investigated the following heat equation
\[
u_t-\triangle u={}_0^{}I_t^{1-\gamma}(|u|^{p-1}u)
\]
on $\mathbb{R}^N$. They obtained the critical exponent of this problem is $\max\{\gamma^{-1},1+\frac{2(2-\gamma)}{N-2(1-\gamma)}\}$.

For problem \eqref{20.1} with $\sigma=1$ and $\mu=0$, Chen and Palmieri \cite{chen} obtained a generalized Strauss exponent $p_0(N,\gamma)$ which is the positive root of the quadratic equation $(N-1)/2p^2-((N+1)/2+1-\gamma)p-1=0$ for $N\geq 2$. Moreover, $p_0(1,\gamma)=\infty$. Under the assumptions that $u_0$ and $u_1$ are nonnegative functions and have compact support,  they proved that the energy solution blows up in finite time if $1<p\leq p_0(N,\gamma)$ for $N\geq 2$, and $p>1$ for $N=1$.

In \cite{a2}, D'Abbicco considered the optimality of the critical exponent for problem \eqref{20.1} with $\sigma=1$, $\eta=\frac{1}{2}$ and $\mu=2$. He proved that if $u_1\geq 0$ is nontrivial and one of the following conditions is satisfied: (a) $N=1$, $\gamma\in (0,1)$ and $p>1$; (b) $N\geq 2$, $\gamma\in [(N-2)/N,1)$ and $1<p\leq 1+(3-\gamma)/(N+\gamma-2)$; (c) $N\geq 3$, $\gamma\in (0,(N-2)/N)$ and $1<p<\gamma^{-1}$, then there is no global weak solution. But the case $\gamma\in (0,(N-2)/N)$ and $p=\gamma^{-1}$ is not included in the above results.
To derive a satisfying nonexistence result for the case $p=\gamma^{-1}$ and $\gamma\in (0,(N-2)/N)$ was left as an open problem (see Remark 4.2 in \cite{a2}).

In \cite{dao}, Dao and Fino studied the nonexistence of global solutions for problem \eqref{20.1} with $\sigma=1$ and $\mu>0$. They proved that if $(u_0,u_1)\in H^{2\eta}(\mathbb{R}^N)\times L^2(\mathbb{R}^N)$, $\int_{\mathbb{R}^N} [u_1(x)+\mu (-\triangle)^\eta u_0(x)]dx>0$, and one of the following conditions is satisfied: (a) $N\leq [2-\gamma(2-\tilde{\eta})]$, $\gamma\in (0,1)$ and $p>1$; (b)$1<p\leq 1+[2+(1-\gamma)(2-\tilde{\eta})]/[N-2+\gamma(2-\tilde{\eta})]_+$ if $[2-\gamma(2-\tilde{\eta})]<N\leq 2$, $\gamma\in (0,1)$ or $N>2$, $\gamma\in [(N-2)/N,1)$; (c) $N\geq 3$, $\gamma\in (0,(N-2)/N)$ and $1<p<\gamma^{-1}$, where $\tilde{\eta}=\min\{2\eta,1\}$, then there exists no global weak solution to this problem. Similarly, the case $\gamma\in (0,(N-2)/N)$ and $p=\gamma^{-1}$ is not included in the their results.

Motivated by the above results, in this paper, we study the nonexistence of global weak solution for problem \eqref{20.1} in the general cases.  We will prove that for problem \eqref{20.1}, the critical case $p\gamma= 1$ is in the nonexistence category.
In particular, we give an answer to an open problem posed in \cite{a2}, and our results complement those obtained in \cite{dao}. Moreover, it should be emphasized that comparing with the above results, any positivity condition is not imposed on the initial values (see Theorem \ref{th1}).

In order to present our results, we need to recall some properties of a test
function proved in \cite{af}.
Let
\begin{equation}\label{11.20}\chi(x)=\Big(\int_{\mathbb{R}^N} (1+|x|^2)^{-\frac{q}{2}}dx\Big)^{-1}(1+|x|^2)^{-\frac{q}{2}},
\end{equation}
where $q$ satisfies that $q>N$ if $\sigma=1$ and $\eta=1$ or $\eta=0$, $q\in (N,N+2\sigma]$ if $0<\sigma<1$ and $\eta=1$ or $\eta=0$, and $q\in (N,N+2\min\{\sigma,\eta\}]$ if $0<\sigma<1$ and $0<\eta<1$. Then it follows from Corollary 3.1 in \cite{af} that
\begin{equation}\label{11.19-2}
|(-\triangle)^\sigma \chi|\leq C\chi(x),\ |(-\triangle)^\eta \chi|\leq C\chi(x)
\end{equation}
for some constant $C>0$.  We will prove the following nonexistence result.
\begin{theorem}\label{th1}
Let $p>1$, $\gamma<1$, $0<\sigma\leq 1$, $0\leq \eta\leq 1$, $\mu\in\mathbb{R}$, and $u_0\chi,u_1\chi\in L^1(\mathbb{R}^N)$. Assume $u$ is a nontrivial weak solution of \eqref{20.1} and $|u|^p\chi \in L_{loc}^1([0,T),L^1(\mathbb{R}^N))$. If $p\gamma\leq 1$, then $T<+\infty$.
\end{theorem}

The rest of the paper is organized as follows.  In Section 2, we present some properties of the
fractional integrals, and study the asymptotic properties of solutions for an integral inequality. In Section 3, we give the proof of Theorem \ref{th1}.

For simplicity, in this paper, we denote by $C$ a positive constant which may differ from one occurrence to the next.
\section{Preliminaries}

\noindent

To prove our results, some preliminaries are presented  in
this section.
For $T>0$, $\alpha>0$, the Riemann-Liouville fractional integrals are defined
by
\[
{_0I_t^{\alpha}}u=\frac{1}{\Gamma(\alpha)}\int_0^t\frac{u(s)}{(t-s)^{1-\alpha}}ds,\ \ {_tI_T^{\alpha}}u=\frac{1}{\Gamma(\alpha)}\int_t^T\frac{u(s)}{(s-t)^{1-\alpha}}ds,
\]
which have the following properties.
\begin{lemma}\label{l1}\cite{kilbas2006,s} Let $\alpha,\beta>0$ and $T>0$.
\begin{enumerate}
  \item [\rm (i)] ${}_0I_t^\alpha$ and ${}_tI_T^\alpha$ are bounded in $L^p(0,T)$ for every $p\geq 1$.
  \item [\rm (ii)] If $f\in L^1(0,T)$, then ${}_0I_t^\alpha({}_0I_t^\beta f)={}_0I_t^{\alpha+\beta}f$ and ${}_tI_T^\alpha({}_tI_T^\beta f)={}_tI_T^{\alpha+\beta}f$.
  \item [\rm (iii)] If $f\in L^p(0,T)$, $g\in L^q(0,T)$ and $p,q\geq1$, $1/p+1/q=1$, then
\begin{equation}\label{112-4}
\int_0^T({_0I_t^\alpha}f)g(t)
dt=\int_0^T({_tI_T^\alpha}g)f(t)dt.
\end{equation}
\item [\rm (iv)] If $f\in L^1(0,\infty)$, then $\mathcal {L}({}_0I_t^\alpha f)=\lambda^{-\alpha}\mathcal {L}f$ for $\text{Re}\lambda>0$ where $\mathcal {L}f$ denotes the Laplace transform of $f$.
\end{enumerate}
\end{lemma}

In order to prove our nonexistence results, drawing on the ideas of proving the blow-up results in  \cite{T.Cazenave}, we study the asymptotic properties of solutions for the following integral inequality.
\begin{lemma}\label{lemma17}
Let $A,B,b\in \mathbb{R}$, $T>0$, $a,c>0$,  $\gamma<1$, $p>1$ and $b^2-4c>0$. If $w\in L_{loc}^p([0,T))$  satisfies that $m(\{t\ |\ w(t)\not=0\})\not =0$ and for almost everywhere $t\in (0,T)$,
\begin{equation}\label{1.7-1}
w+b\cdot\prescript{}{0}{I}_t^{1}w+c\cdot\prescript{}{0}{I}_t^{2}w-A-Bt\geq a\cdot\prescript{}{0}{I}_t^{3-\gamma}|w(t)|^p,
\end{equation}
then the following properties hold.
\begin{enumerate}
\item[\rm (i)] For $l>\frac{p(3-\gamma)}{p-1}$, there exists a positive constant $K$ such that \[A(l+\gamma-1)+BT\leq K(T^{-\frac{3-\gamma}{p-1}}+ T^{1-\frac{2-\gamma}{p-1}}+T^{2-\frac{1-\gamma}{p-1}}),\]
      where $K$ depends on $l$ but not on $T$.

\item[\rm(ii)] If $T=+\infty,$ then $\liminf_{t\rightarrow +\infty}t^{\gamma-2} \prescript{}{0}{I}_t^{2}w>0$.

\item[\rm(iii)] If $p\gamma\leq 1$, then $T<+\infty$.
\end{enumerate}
\end{lemma}
\begin{proof} (i) According to \eqref{112-4} and \eqref{1.7-1}, we know that
\begin{equation}\label{6.10-2}
\int_0^Tw[\varphi +b\cdot\prescript{}{t}{I}_T^{1}\varphi+c\cdot\prescript{}{t}{I}_T^{2}\varphi]dt\geq a\int_0^T|w|^p\prescript{}{t}{I}_T^{3-\gamma}\varphi dt+A\int_0^T\varphi dt+B\int_0^Tt\varphi dt,
\end{equation}
for every $\varphi\in C([0,T])$ with $\varphi\geq0$.
It is well-known that for $l>\frac{p(3-\gamma)}{p-1}$ and $\beta>0$
\begin{equation}\label{113-3}
\frac{\Gamma(l+1)}{\Gamma(l+\gamma-2)} \prescript{}{t}{I}_T^{\beta}(1-\frac{t}{T})^{l-(3-\gamma)}=\frac{\Gamma(l+1)}{\Gamma(l+\gamma+\beta-2)}
T^{\beta}
(1-\frac{t}{T})^{l-(3-\gamma)+\beta},
\end{equation}
(see, e.g., \cite{kilbas2006,s}).
We denote $\psi_T(t)=(1-\frac{t}{T})^l$ and take $\varphi(t)=\frac{\Gamma(l+1)}{\Gamma(l+\gamma-2)}T^{-(3-\gamma)}\psi_T^\frac{l-(3-\gamma)}{l}$, and then we deduce from \eqref{6.10-2}-\eqref{113-3} and Young's inequality that there exists a constant $C>0$ such that
\begin{align*}
&A\int_0^T\varphi dt
+B\int_0^Tt\cdot \varphi dt +a\int_0^T|w|^p\psi_Tdt
\leq \int_0^Tw[\varphi
+b\cdot\prescript{}{t}{I}_T^{1}\varphi
+c\cdot\prescript{}{t}{I}_T^{2}\varphi]dt\\
&\leq \frac{a}{2}\int_0^T|w|^p \psi_T dt
+CT^{1-\frac{p(3-\gamma)}{p-1}}+ CT^{1-\frac{p(2-\gamma)}{p-1}}+CT^{1-\frac{p(1-\gamma)}{p-1}}.
\end{align*}
This implies that
\begin{align}\label{11.1-1}
&\frac{a}{2}\int_0^T|w|^p \psi_T dt+A\frac{\Gamma(l+1)}{\Gamma(l+\gamma-1)}T^{\gamma-2}
+B\frac{\Gamma(l+1)}{\Gamma(l+\gamma)}T^{\gamma-1}\nonumber\\
&\leq CT^{1-\frac{p(3-\gamma)}{p-1}}+ CT^{1-\frac{p(2-\gamma)}{p-1}}+CT^{1-\frac{p(1-\gamma)}{p-1}}.
\end{align}
As a result, there exists a constant $K>0$ such that
\begin{align*}
A(l+\gamma-1)+BT\leq K(T^{-\frac{3-\gamma}{p-1}}+ T^{1-\frac{2-\gamma}{p-1}}+T^{2-\frac{1-\gamma}{p-1}}).
\end{align*}

(ii) For $\tilde{T}<+\infty$ and $t\in [0,\tilde{T}]$, we denote $f(t)=w+b\cdot\prescript{}{0}{I}_t^{1}w+c\cdot\prescript{}{0}{I}_t^{2}w-A-Bt$. We set $f(t)=0$ for $t>\tilde{T}$. Then $f\in L^1(0,+\infty)$. Consider the integral equation
\[w+b\cdot\prescript{}{0}{I}_t^{1}w+c\cdot\prescript{}{0}{I}_t^{2}w-A-Bt=f(t), \ t\in [0,+\infty).\]
Setting  $v(t)=\prescript{}{0}{I}_t^2 w$, then using Lemma \ref{l1} we know that $v(t)$ satisfies
\begin{equation}\label{11.9-1}
v+b\cdot\prescript{}{0}{I}_t^{1}v+c\cdot\prescript{}{0}{I}_t^{2}v-\frac{At^2}{2}-\frac{Bt^3}{6}
=\prescript{}{0}{I}_t^{2}f.
\end{equation}
Applying the Laplace transform to \eqref{11.9-1}, we get from Lemma \ref{l1} that for $\text{Re}s>0$
\[
\mathcal {L} v+\frac{b}{s}\mathcal {L} v+\frac{c}{s^2}\mathcal {L} v-\frac{1}{s^2}(\frac{A}{s}+\frac{B}{s^2})=\frac{1}{s^2}\mathcal {L}f,
\]
that is $\mathcal {L} v=\frac{1}{s^2+bs+c}(\frac{A}{s}+\frac{B}{s^2})+\frac{\mathcal {L}f}{s^2+bs+c}$. Since $\frac{1}{\lambda_2-\lambda_1}\mathcal {L} (e^{\lambda_2t}-e^{\lambda_1t})=\frac{1}{s^2+bs+c}$, where
$\lambda_1=\frac{-b-\sqrt{b^2-4c}}{2}$ and $\lambda_2=\frac{-b+\sqrt{b^2-4c}}{2}$,
then we deduce from the uniqueness of Laplace transform  that for $t\in [0,+\infty)$
\begin{equation}\label{11.21-1}
v=\frac{1}{\lambda_2-\lambda_1}\big[A\cdot\prescript{}{0}{I}_t^{1}(e^{\lambda_2t}-e^{\lambda_1t})
+B\cdot\prescript{}{0}{I}_t^{2}(e^{\lambda_2t}-e^{\lambda_1t})
+\int_0^t(e^{\lambda_2(t-s)}-e^{\lambda_1(t-s)})f(s)ds\big].
\end{equation}
Since $f(t)\geq a\cdot\prescript{}{0}{I}_t^{3-\gamma}|w(t)|^p$ for $t\in (0,\tilde{T})$ and $e^{\lambda_2t}-e^{\lambda_1t}\geq 0$, by the arbitrariness of $\tilde{T}$, we deduce from  \eqref{11.21-1} and
Fubini's theorem that for $t\geq 0$
\begin{align}\label{11.9-2}
\prescript{}{0}{I}_t^2 w=v(t)
\geq &\frac{A}{\lambda_2-\lambda_1}\prescript{}{0}{I}_t^{1}(e^{\lambda_2t}-e^{\lambda_1t})
+\frac{B}{\lambda_2-\lambda_1}\prescript{}{0}{I}_t^{2}(e^{\lambda_2t}-e^{\lambda_1t})\nonumber\\
&+\frac{a}{\lambda_2-\lambda_1}\int_0^t(e^{\lambda_2(t-s)}-e^{\lambda_1(t-s)})
\prescript{}{0}{I}_s^{3-\gamma}|w(s)|^pds\nonumber\\
=&\frac{A}{\lambda_2-\lambda_1}\prescript{}{0}{I}_t^{1}(e^{\lambda_2t}-e^{\lambda_1t})
+\frac{B}{\lambda_2-\lambda_1}\prescript{}{0}{I}_t^{2}(e^{\lambda_2t}-e^{\lambda_1t})\nonumber\\
&+\frac{a}{\lambda_2-\lambda_1}\int_0^t\prescript{}{0}{I}_{t-s}^{3-\gamma}[e^{\lambda_2(\cdot)}-e^{\lambda_1(\cdot)}]
|w(s)|^pds
\end{align}

In terms of our assumption that $m(\{t\ |\ w(t)\not=0\})\not =0$, there exist $t_0\geq0$ and $\delta>0$ such that $\int_{t_0}^{t_0+\delta}|w(t)|^pdt\not=0$. On the other hand, by some simple calculations, we see that for $\lambda <0$
\begin{equation}\label{11.9-4}
\lim_{t\rightarrow +\infty} \prescript{}{0}{I}_t^{1}e^{\lambda t}=\lim_{t\rightarrow +\infty} t^{-1} \prescript{}{0}{I}_t^{2}e^{\lambda t}=-\frac{1}{\lambda}, \lim_{t\rightarrow +\infty} t^{\gamma-2} \prescript{}{0}{I}_t^{3-\gamma}e^{\lambda t}=-\frac{1}{\lambda\Gamma(3-\gamma)}.
\end{equation}
Then it follows from \eqref{11.9-2} that for $t$ large enough,
\begin{align}\label{11.9-3}
\prescript{}{0}{I}_t^2 w(t)\geq& \frac{A}{\lambda_2-\lambda_1}\prescript{}{0}{I}_t^{1}(e^{\lambda_2t}-e^{\lambda_1t})
+\frac{B}{\lambda_2-\lambda_1}\prescript{}{0}{I}_t^{2}(e^{\lambda_2t}-e^{\lambda_1t})\nonumber\\
&+ \frac{a}{\lambda_2-\lambda_1} \int_{t_0}^{t_0+\delta}\prescript{}{0}{I}_{t-s}^{3-\gamma}[e^{\lambda_2(\cdot)}-e^{\lambda_1(\cdot)}]|w(s)|^pds\nonumber\\
\geq&\frac{A}{\lambda_2-\lambda_1}\prescript{}{0}{I}_t^{1}(e^{\lambda_2t}-e^{\lambda_1t})
+\frac{B}{\lambda_2-\lambda_1}\prescript{}{0}{I}_t^{2}(e^{\lambda_2t}-e^{\lambda_1t})\nonumber\\
&+\frac{a}{\lambda_2-\lambda_1}\prescript{}{0}{I}_{t-t_0-\delta}^{3-\gamma}
[e^{\lambda_2(\cdot)}-e^{\lambda_1(\cdot)}]\int_{t_0}^{t_0+\delta}|w(s)|^pds\nonumber\\
\geq& \frac{A}{\lambda_2-\lambda_1}\prescript{}{0}{I}_t^{1}(e^{\lambda_2t}-e^{\lambda_1t})
+\frac{B}{\lambda_2-\lambda_1}\prescript{}{0}{I}_t^{2}(e^{\lambda_2t}-e^{\lambda_1t})
+\frac{at^{2-\gamma}}{2c\Gamma(3-\gamma)}\int_{t_0}^{t_0+\delta}|w|^pds,
\end{align}
where we have used the fact that $\prescript{}{0}{I}_{t}^{3-\gamma}
[e^{\lambda_2(\cdot)}-e^{\lambda_1(\cdot)}]$ is an increasing function on $[0,+\infty)$.
Recalling $\gamma<1$, we then derive from \eqref{11.9-4} and \eqref{11.9-3} that $\liminf_{t\rightarrow +\infty}t^{\gamma-2} \prescript{}{0}{I}_t^{2}w>0$.

(iii) If the assertion would not hold, then it follows from \eqref{11.1-1} that for every $T>0$,
\begin{equation}\label{11.12-2}
\int_0^T|w|^p \psi_T dt\leq C[T^{1-\frac{p(3-\gamma)}{p-1}}+T^{1-\frac{p(2-\gamma)}{p-1}}+T^{1-\frac{p(1-\gamma)}{p-1}}
+|A|T^{\gamma-2}+|B|T^{\gamma-1}].
\end{equation}

If $p\gamma<1$, then $1-\frac{p(1-\gamma)}{p-1}<0$. Hence using \eqref{11.12-2} and taking $T\rightarrow +\infty$, we see that  $\int_0^{+\infty}|w|^p  dt=0$. This implies $w\equiv0$. We obtain a contradiction.

If $p\gamma=1$, then by letting $T\rightarrow +\infty$, it follows from \eqref{11.12-2} that $\int_0^{+\infty}|w|^p  dt<+\infty.$ On the other hand, Property (ii) yields that there exist constants $C>0$ and $t_1>0$ such that $\prescript{}{0}{I}_t^{2}w\geq Ct^{2-\gamma}$ for $t\geq t_1$.
Then using the following weighted estimate of the operator $\prescript{}{0}{I}_t^{2}$
\[
\Big(\int_0^{+\infty}t^{-2 p} |\prescript{}{0}{I}_t^{2} w|^p dt\Big)^{\frac{1}{p}} \leq \frac{\Gamma(1-1/ p)}{\Gamma(3-1/p)}\Big(\int_0^{+\infty}|w(t)|^p  dt\Big)^{\frac{1}{p}}
\]
(see, e.g., Theorem
3.7 in \cite{s}), we have
\[
C\int_{t_1}^{+\infty}t^{-1} dt= C\int_{t_1}^{+\infty}t^{-2 p+p(2-\gamma)} dt\leq\Big(\frac{\Gamma(1-1/ p)}{\Gamma(3-1/p)}\Big)^p\int_0^{+\infty}|w(t)|^p  dt<+\infty.
\]
We again get a contradiction. Hence $T<+\infty$.
\end{proof}

Finally, we give the definition of  weak solution to problem \eqref{20.1}.
\begin{definition}\label{de1}
For $T>0$, $(1+|x|)^{-N-2\eta}u_0\in L^{1}(\mathbb{R}^N)$ and $u_1\in L_{loc}^1(\mathbb{R}^N)$, we call that $u$ is a weak solution of \eqref{20.1} if $u\in L^p((0,T),L_{loc}^p(\mathbb{R}^N))$, $(1+|x|)^{-N-2\min\{\eta,\sigma\}}u\in L^{1}((0,T),L^{1}(\mathbb{R}^N))$ and
\begin{align*}
&\int_{\mathbb{R}^N}\int_0^T[{_0I_t^{1-\gamma}}(|u|^{p})\varphi+(u_0+tu_1 )\varphi_{tt}+\mu u_0 (-\triangle)^\eta(-\varphi_t)]dtdx\\
&=\int_{\mathbb{R}^N}\int_0^Tu(-\triangle)^\sigma\varphi dtdx
+\int_{\mathbb{R}^N}\int_0^Tu[\varphi_{tt}+\mu (-\triangle)^\eta(-\varphi_t)] dtdx
\end{align*}
for any test function $\varphi(t,x)=\psi(x)\phi(t)$ with $\psi\in C_0^\infty (\mathbb{R}^N)$, $\phi\in C^2([0,T])$ and $\phi(T)=\phi'(T)=0$. Furthermore, we say that $u$ is a global weak solution of \eqref{20.1} if $T > 0$ can be arbitrarily chosen.
\end{definition}
\begin{remark}
In terms of Proposition 2.9 in \cite{g}, we know that for some constant $C>0$,
\[|(-\triangle)^\sigma \psi(x)|\leq C(1+|x|^2)^{-\frac{N+2\sigma}{2}},\ |(-\triangle)^\eta \psi(x)|\leq C(1+|x|^2)^{-\frac{N+2\eta}{2}}.\]
Hence the integrals in Definition \ref{de1} are well-defined.
\end{remark}

\section{The proof of Theorem \ref{th1}}
\noindent

In this section, we will prove Theorem \ref{th1}.

\noindent\textbf{Proof of Theorem \ref{th1}} We only give the proof for the case $0<\sigma,\eta<1$, since if one of $\sigma$ and $\eta$ is an integer, the conclusion can be derived by some similar arguments .
Suppose that $u$ is  a global weak
solution of \eqref{20.1}. Take $\Psi\in C_0^\infty(\mathbb{R}^N)$ such that
$\Psi(x)=1$ for $|x|<1$ and $\Psi(x)=0$ for $|x|>2$. Let $T>0$ and $\psi_n(x)=\Psi(\frac{x}{n})$, $n=1,2,\dots$.
Choosing the test function $\varphi=\psi_n(x)\chi(x)\phi(t)$ in Definition \ref{de1}, where $\chi$ is given by \eqref{11.20} and $\phi\in C^2([0,T])$ with $\phi(t)\geq 0$ and $\phi(T)=\phi'(T)=0$, we obtain that
\begin{align}\label{11.17-1}
&\int_{\mathbb{R}^N}\int_0^T[{_0I_t^{1-\gamma}}(|u|^{p})\psi_n\chi\phi+(u_0+tu_1 )\psi_n\chi\phi''+\mu u_0 (-\triangle)^\eta(\psi_n\chi)(-\phi')]dtdx\nonumber\\
&=\int_{\mathbb{R}^N}\int_0^Tu(-\triangle)^\sigma(\psi_n\chi)\phi dtdx
+\int_{\mathbb{R}^N}\int_0^Tu[\psi_n\chi\phi''+\mu (-\triangle)^\eta(\psi_n\chi)(-\phi') ]dtdx
\end{align}

Since for $s\in (0,1)$,
\begin{equation}\label{11.25-1}
\Big|\frac{2(\psi_n\chi)(x)-(\psi_n\chi)(x+y)-(\psi_n\chi)(x-y)}{|y|^{N+2s}}\Big|
\leq \left\{\begin{array}{l}\frac{C}{|y|^{N+2s-2}},\
|y|\leq 1,\\
\frac{C}{|y|^{N+2s}},\ |y|>1,
\end{array}\right.
\end{equation}
we deduce from the Lebesgue dominant convergence theorem that
\[
(-\triangle)^s(\psi_n\chi)=\frac{s2^{2s-1}\Gamma(\frac{N+2s}{2})}{\pi^{\frac{N}{2}}\Gamma(1-s)}
\int_{\mathbb{R}^N} \frac{2(\psi_n\chi)(x)-(\psi_n\chi)(x+y)-(\psi_n\chi)(x-y)}{|y|^{N+2s}} dy\rightarrow (-\triangle)^s \chi.
\]

The proof of Proposition 2.9 in \cite{g} yields that for $s\in (0,1)$ and $|x|>1$, $|(-\triangle)^s(\psi_n\chi)(x)|\leq C_{N,s,\phi_n\chi}|x|^{-N-2s}$, where
\begin{equation}\label{11.24-1}
C_{N,s,\phi_n\chi}=\Big(\sup_{x\in \mathbb{R}^N}(1+|x|)^{N+2}\sum_{|\alpha|=2}|\partial^\alpha (\psi_n\chi)|
+\sup_{x\in \mathbb{R}^N}(1+|x|)^{N}|\psi_n\chi|+\|\psi_n\chi\|_{L^1(\mathbb{R}^N)}
\Big).
\end{equation}

Next, we prove that $|(-\triangle)^s(\psi_n\chi)(x)|\leq C(1+|x|)^{-N-2s}$ for some positive constant $C$ independent of $n$.  In fact, since $q>N$,
\[|(\phi_n)_{x_ix_j}\chi|\leq \frac{C}{n^2}\chi_{\{x|n<|x|<2n\}}(x)(1+|x|^2)^{-q/2}\leq C(1+|x|^2)^{-q/2-1},\]
where $\chi_A(x)$ denotes the characteristic function of the set $A$,
$|\chi_{x_ix_j}\phi_n|\leq C(1+|x|^2)^{-q/2-1}$ and
\[|(\phi_n)_{x_j}\chi_{x_i}+(\phi_n)_{x_i}\chi_{x_j}|\leq \frac{C}{n}\chi_{\{x|n<|x|<2n\}}(x)(1+|x|^2)^{-(q+1)/2}\leq C(1+|x|^2)^{-q/2-1},\]
we know that $\sup_{x\in \mathbb{R}^N}(1+|x|)^{N+2}\sum_{|\alpha|=2}|\partial^\alpha (\psi_n\chi)|\leq C$.
On the other hand, we can easily see that $\sup_{x\in \mathbb{R}^N}(1+|x|)^{N}|\psi_n\chi|+\|\psi_n\chi\|_{L^1(\mathbb{R}^N)}\leq C$. Thus $C_{N,s,\phi_n\chi}\leq C$.
Furthermore, by \eqref{11.25-1}, we can obtain the desired conclusion.

Letting $n\rightarrow +\infty$ and using \eqref{11.17-1}, \eqref{11.25-1} and \eqref{11.19-2}, we deduce from the Fatou Lemma and the Lebesgue dominant convergence theorem that
\begin{align}\label{11.19-1}
&\int_{\mathbb{R}^N}\int_0^T[{_0I_t^{1-\gamma}}(|u|^{p})\chi\phi+(u_0+tu_1 )\chi\phi''+\mu u_0 (-\triangle)^\eta(\chi)(-\phi')]dtdx\nonumber\\
&\leq \int_{\mathbb{R}^N}\int_0^Tu(-\triangle)^\sigma(\chi)\phi dtdx
+\int_{\mathbb{R}^N}\int_0^Tu[\chi\phi''+\mu (-\triangle)^\eta(\chi)(-\phi') ]dtdx\nonumber\\
&\leq C\int_{\mathbb{R}^N}\int_0^T|u|\chi\phi dtdx
+\int_{\mathbb{R}^N}\int_0^T\chi|u|\cdot|\phi''|dtdx+|\mu| C \int_{\mathbb{R}^N}\int_0^T\chi|u|\cdot|-\phi'|dtdx.
\end{align}
Denote $w(t)=\int_{\mathbb{R}^N} |u|\chi dx$. Our assumptions yield $w\in L^p(0,T)$. Then it follows from \eqref{11.19-1} and Jensen's inequality that
\begin{align}\label{11.19-5}
&\int_0^T{_0I_t^{1-\gamma}}(w^{p})\phi dt+\int_{\mathbb{R}^N}\int_0^T[(u_0+tu_1 )\chi\phi''+\mu u_0 (-\triangle)^\eta(\chi)(-\phi')]dtdx\nonumber\\
&\leq C\int_0^Tw\phi dt
+\int_0^Tw|\phi''|dtdx+|\mu| C \int_0^Tw|-\phi'|dtdx
\end{align}
Taking $\phi(t)={}_tI_T^{2}\tilde{\phi}$ with $\tilde{\phi}\in C_0^\infty((0,T))$ and $\tilde{\phi}\geq 0$ in \eqref{11.19-5}, we deduce from Lemma \ref{l1} that
\begin{align*}
&\int_0^T{}_0I_t^{3-\gamma}w^{p}\tilde{\phi}dt+\int_{\mathbb{R}^N}\int_0^T[(u_0+tu_1 )\chi\tilde{\phi}+\mu t u_0 (-\triangle)^\eta(\chi)\tilde{\phi}]dtdx\\
&\leq C\int_0^T{}_0I_t^{2}w\tilde{\phi} dt
+\int_0^Tw\tilde{\phi} dt+C \int_0^T{}_0I_t^{1}w\tilde{\phi}dt,
\end{align*}
where we have used the fact that $\frac{d}{dt}({}_tI_T^{2}\tilde{\phi})=-{}_tI_T^{1}\tilde{\phi}$ and $\frac{d^2}{dt^2}({}_tI_T^{2}\tilde{\phi})=\tilde{\phi}$.
Noting that $w\geq 0$, it follows that
\[
\int_0^T[{}_0I_t^{3-\gamma}w^{p}+A+Bt]\tilde{\phi}dt\leq  C\int_0^T{}_0I_t^{2}w\tilde{\phi} dt+\int_0^Tw\tilde{\phi} dt+(2\sqrt{C}+1) \int_0^T{}_0I_t^{1}w\tilde{\phi}dt
\]
where $A=\int_{\mathbb{R}^N}u_0\chi dx$, $B=\int_{\mathbb{R}^N}[u_1\chi+\mu u_0(-\triangle)^\eta\chi]dx$. Due to the arbitrariness of $\tilde{\phi}$, we have
\[
{}_0I_t^{3-\gamma}w^{p}\leq w+(2\sqrt{C}+1){}_0I_t^{1}w+C({}_0I_t^{2}w)-A-Bt.
\]
Therefore we get a contradiction by Lemma \ref{lemma17}(iii).

\section*{Acknowledgment}
\noindent

This work was supported in part by Young Backbone Teachers of Henan Province (No.
2021GGJS130).


\begin{thebibliography}{99}

\bibitem{Berbiche2}M. Berbiche, Existence and blow-up of solutions for damped wave system with nonlinear
memory, Appl. Anal. 94(12) (2015), 2535-2564.

\bibitem{T.Cazenave} T. Cazenave, F. Dickstein, F.B. Weissler, An
equation whose Fujita critical exponent is not given by scaling,
Nonlinear Anal. 68 (2008), 862-874.

\bibitem{chen} W. Chen, Interplay effcts on blow-up of weakly coupled systems for
semilinear wave equations with general nonlinear memory terms, Nonlinear Anal. 202 (2021), Article 112160.

\bibitem{chen2} W. Chen, A. Palmieri, in: M. Cicognani, D. Del Santo, A. Parmeggiani, M. Reissig (Eds.), Blow-up Result for a Semilinear Wave Equation with a Nonlinear Memory Term, in: Springer INdAM Series, vol. 43, 2020, p. 20.

\bibitem{a1} M. D'Abbicco, The influence of a nonlinear memory on the damped
wave equation, Nonlinear Anal. 95 (2014), 130-145.

\bibitem{a2} M. D'Abbicco, A wave equation with structural damping
and nonlinear memory, Nonlinear Differ. Equ. Appl. 21 (2014), 751-773.

\bibitem{af} M. D'Abbicco, K. Fujiwara,
A test function method for evolution equations with fractional
powers of the Laplace operator, Nonlinear Anal. 202 (2021), Article 112114.


\bibitem{Andrade} B. de Andrade, A. Viana, Integrodifferential equations with applications to a plate equation with memory, Math. Nachr. 289 (17-18) (2016), 2159-2172.

\bibitem{dao} T.A. Dao, A.Z. Fino,  Blow-up results for a semi-linear structural damped wave model with nonlinear memory, Math. Nachr. 295 (2)(2022), 309-322.

\bibitem{fino}A.Z. Fino, Critical exponent for damped wave equations with nonlinear memory. Nonlinear Anal. 74(16) (2011), 5495-5505.

\bibitem{g} N. Garofalo, Fractional thoughts. New developments in the analysis of nonlocal operators, Contemp. Math. 723, 1-135, Amer. Math. Soc., Providence, RI, 2019.

\bibitem{kilbas2006} A.A. Kilbas, H.M. Srivastava, J.J. Trujillo, Theory and
Applications of Fractional Differential Equations, vol 204. Elsevier
Science B.V., Amsterdam, 2006.

\bibitem{s} S.G. Samko, A.A. Kilbas, O.I. Marichev, Fractional Integrals and Derivatives:
Theory and Applications, Gordon and Breach, Yverdon, 1993.

\end{thebibliography}
\end{document}